\documentclass[12pt]{article}
\usepackage{mathpple}
\usepackage{charter}
\usepackage{amsfonts,amssymb, lineno}
\usepackage{amsmath}
\usepackage{amssymb}
\usepackage{amsfonts}
\usepackage{graphicx}
\usepackage{authblk} 
\usepackage{biblatex}

\newtheorem{corollary}{Corollary}

\newtheorem{lemma}{Lemma}

\newtheorem{theorem}{Theorem}
\numberwithin{equation}{section}

\begin{document}
\title{On the locating-chromatic number of corona product of graphs}
\author [1] { D.K. Syofyan } 
\author [2] { S.W. Saputro } 
\author [2] { E.T. Baskoro }
\author [1] { I.A. Purwasih }
\affil [1] { \textit {Combinatorial Mathematics Research Group}\\
\textit {Computing Research Center}\\
\textit {National Research and Innovation Agency}\\
\textit {Jalan Raya Jakarta-Bogor KM 46, Cibinong, Bogor, Jawa Barat, Indonesia, 16911}\\
\textit {Email: dian.kastika.syofyan@brin.go.id}} 
\affil [2] { \textit {Combinatorial Mathematics Research Group}\\
\textit {Faculty of Mathematics and Natural Sciences}\\
\textit {Institut Teknologi Bandung}\\
\textit {Jalan Ganesa 10 Bandung 40132, Indonesia.} } 
\date{}
\maketitle

\begin{abstract}
Let $G=(V,E)$ be a finite, simple, and connected graph. The \textit{locating-chromatic number} of a graph $G$ can be defined as the cardinality of a minimum resolving partition of the vertex set $V(G)$ such that all vertices have different coordinates and every two adjacent vertices in $G$ is not contained in the same partition class. In this case, the coordinate of a vertex in $G$ is expressed in terms of the distances of this vertex to all partition classes. The \textit{corona product} of a graph $G$ of order $n$ and a graph $H,$ denoted by $G \odot H,$ is the graph obtained by taking one copy of $G$ and $n$ copies of $H$ and joining the $i^{th}$-vertex of $G$ to every vertex in the $i^{th}$-copy of $H$. In this paper, we determine the sharp general bound of the locating-chromatic number of $G \odot H$ for $G$ is a connected graph and $H$ is an arbitrary graph, or $G$ is a tree graph and $H$ is a complement of complete graph.
\end{abstract}

\maketitle



\section{Introduction}

The concept of locating-chromatic number was introduced by Chartrand \emph{et al.} in 2002 \cite{chartrand2002} as a special case of partition dimension concept \cite{Chartrand1998}. They have provided the boundaries for the the locating-chromatic number of a graph which depends on order and diameter of its graph. They also have determined the locating-chromatic number of paths, cycles, and double stars.

Throughout this paper, all graphs $G$ are finite and simple. \ Let $G=(V,E)$ be a connected graph. The \emph{ $k$-coloring} of $G$ is a function $c:V(G)\rightarrow \{1,2,\ldots,k\}$ where $c(u)\neq c(v)$ for two adjacent vertices $u$ and $v$ in $G$. Let $\Pi = \{C_1, C_2, \ldots , C_k\}$ be the partition induced by a $k$-coloring $c$ on $V(G)$. For $1\leq i\leq k$, $C_i$ is called a \emph{color class} of $G$. The \emph{color code} $c_{\Pi} (v)$ of a vertex $v$ in $G$ is defined as $(d(v, C_1), d(v, C_2), \ldots, d(v, C_k))$ where $d(v, C_i)= \min\{d(v,x) | x \in  C_i\}$ for $ 1 \leq  i \leq  k$. If all distinct vertices of $G$ have distinct color codes, then $c$ is called a \emph{locating coloring} of $G$. The least integer $k$ such that there is a locating coloring in $G$ is called the \emph{locating-chromatic number} of $G,$ denoted by $\chi_L(G)$.

The locating-chromatic number of some well-known classes of graph have been obtained by some authors. Asmiati \emph{et al.} \cite{Asmiati2011,Asmiati20121} have determined the locating-chromatic number of amalgamation of stars and firecrackers. In \cite{Behtoei20111}, Behtoei \emph{et al.} studied the locating-chromatic number of Kneser graph. Meanwhile, Purwasih \emph{et al.} \cite{Purwasih20121}  have determined the locating-chromatic number of Halin graph. Some authors also have determine the locating-chromatic number of graphs obtained from a product graphs.  The locating-chromatic number of Cartesian product of graphs and joint product graphs have been proved by Behtoei \emph{et al.} \cite{Behtoei20112,Behtoei20113}. Purwasih \emph{et al.} \cite{Purwasih20122} obtained the locating-chromatic number of strong product graphs.

In \cite{Baskoro2012}, Baskoro and Purwasih have determined the upper bound of the locating-chromatic number of corona product of two connected graphs $G$ and $H,$ where the diameter of $H$ is at most two. In this paper, we generalized the results of Baskoro and Purwasih for a connected graph $G$ and an arbitrary graph $H.$

Let $G$ be a connected graph of order $n$ and $V(G)=\{a_1,a_2,\ldots,a_n\}$. The \emph{corona product} between $G$ and $H$, denoted by $G \odot H$, is the graph obtained by taking one copy of $G$ and $n$ copies of $H$ and joining the vertex $a_i$ of $G$ to every vertex in the $i^{th}$-copy of $H$.  By the definition of corona product, we define the vertex set $V(G\odot H)=U\cup V$ where $U$ and $V$ are vertices of $G\odot H$ from $G$ and $n$ copies of $H$, respectively. Let $U=\{(u)\mid u\in V(G)\}$ and $H(u)=\{(u,v)\mid v\in V(H)\}$ be a vertex set from a copy of $H$ whose all vertices are adjacent to $(u)$. So, $V=\bigcup_{u\in V(G)} H(u)$.

Let $H$ be an arbitrary graph containing $k\geq1$ components $H_1,H_2,\ldots,H_k$. For $u\in V(G)$ and $1\leq t\leq k$, we define $H_t(u)=\{(u,v)\mid v\in V(H_t)\}$. So, for every $u\in V(G)$, $H(u)=\bigcup_{1\leq t\leq k} H_t(u)$. We also use some following definitions. Let $c$ be a locating coloring of $G\odot H$ and $\Pi$ be a partition of $V(G\odot H)$ induced by $c$. For every $u\in V(G)$ and $1\leq t\leq k$, let $\Pi(u)$ and $\Pi_t(u)$ be partitions of $H(u)$ and $H_t(u)$, respectively, induced by $c$. Note that, $\Pi_t(u)\subseteq\Pi(u)\subseteq\Pi$.

To prove some results in this paper, we use the following lemma and corollary which are useful to determine the locating-chromatic number of a graph $G$.

\begin{lemma}\label{lema}
Let $G$ be a connected non trivial graph. Let $c$ be a locating coloring of $G$ and $u, v \in V (G)$. If $d(u,w) = d(v,w)$ for every $w \in V (G)\setminus \{u, v\},$ then the color of $u$ and $v$ must be different.
\end{lemma}

\begin{corollary}\label{ColA}
If $G$ is a connected graph containing a vertex which is adjacent to $k$ endpoints of $G$, then $\chi_L(G)\geq k+1$.
\end{corollary}

\section{The Boundaries of Locating-Chromatic Number of $G\odot H$}

For $u\in V(G)$ and $1\leq t\leq k$, let us consider $H_t(u)$ and a vertex $(u)$. By the definition of $G\odot H$, an induced subgraph of $G\odot H$ by $H_t(u)$ and $(u)$ is isomorphic to a joint graph $H_t+K_1$. In lemma below, we prove that $H_t(u)$ is partitioned into at least $\chi_L(H_t+K_1)-1$ color classes induced by a locating coloring of $G\odot H$.

\begin{lemma}\label{join}
Let $G$ be a connected graph of order $n\geq2$ and $H$ be an arbitrary graph containing $k$ components $H_1,H_2,\ldots,H_k$. Let $c$ be a locating coloring of $G\odot H$. For $u\in V(G)$ and $1\leq t\leq k$, the vertex set $H_t(u)$ is partitioned into at least $\chi_L(H_t+K_1)-1$ color classes induced by $c$.
\end{lemma}

\noindent \emph{Proof} :

Let $Q$ be a graph induced by $H_t(u)\cup\{(u)\}$ which is isomorph to $H_t+K_1$. Then $Q$ must be partitioned into at least $\chi_L(H_t+K_1)$ color classes induced by $c$. Since $(u)$ is adjacent to every vertex in $H_t(u)$, the color of $(u)$ must be different than the color of all vertices in $H_t(u)$. Therefore, $H_t(u)$ is partitioned into at least $\chi_L(H_t+K_1)-1$ color classes. 

In two lemmas below, we give a lower bound and the upper bound of $\chi_L(G\odot H)$, respectively.

\begin{lemma}\label{teo1bawah}
Let $G$ be a connected graph of order $n\geq2$ and $H$ be an arbitrary graph containing $k$ components $H_1,H_2,\ldots,H_k$. Then $\chi_L(G\odot H)\geq\max\{\chi_L(H_t+K_1)\mid 1\leq t\leq k\}$.
\end{lemma}

\noindent \emph{Proof} :

Let $c$ be a locating coloring of $G\odot H$. For $t\in\{1,2,\ldots,k\}$ and $u\in V(G)$, by Lemma \ref{join} and considering that $(u)$ is adjacent to every vertex of $H_t(u)$, the vertex set $H_t(u)\cup\{(u)\}$ is partitioned into at least $\chi_L(H_t+K_1)$ color classes induced by $c$. However, it is also possible to have two conditions as follows.
\begin{itemize}
\item For $q\in\{1,2,\ldots,k\}$ and $v\in V(G)$, it is possible to have a color class $C$ in both $H_t(u)$ and $H_q(v)$.
\item For $q\in\{1,2,\ldots,k\}\setminus\{t\}$, if $\chi_L(H_t+K_1)<\chi_L(H_q+K_1)$, then $\chi_L(H_q+K_1)$ cannot be partitioned into $\chi_L(H_t+K_1)$ color classes induced by $c$.
\end{itemize}
Therefore, we obtain that $\chi_L(G\odot H)\geq\max\{\chi_L(H_t+K_1)\mid 1\leq t\leq k\}$. 

\begin{lemma}\label{teo1atas}
Let $G$ be a connected graph of order $n\geq2$ and $H$ be an arbitrary graph containing $k$ components $H_1,H_2,\ldots,H_k$. Then $\chi_L(G\odot H)\leq \chi_L(G)+\sum_{t=1}^k(\chi_L(H_t+K_1)-1)$.
\end{lemma}

\noindent \emph{Proof} :

Let $V(G)=\{u_1,u_2,\ldots,u_n\}$. Let $\chi_L(G)=l$ and $f$ be a locating coloring of $G$ with $l$ colors. For $1\leq t\leq k$, let $\chi_L(H_t+K_1)=m_t$ and $c_t$ be a locating coloring of $H_t+K_1$ with $m_t$ colors such that vertex of $K_1$ colored by $m_t$ and vertices of $H_t$ colored by $1$ until $m_t-1$. Define $c:V(G\odot H)\rightarrow \{1,2,\ldots,l+\sum_{t=1}^{k}(m_t-1)\}$ as follows:

$c((x))=f(x)$ for $x\in V(G)$, and

$c((x,y))=\left\{
\begin{array}{ll}
c_1(y)+l, & \text{for }t=1,\\
c_t(y)+l+\sum_{j=1}^{t-1}(m_j-1), & \text{for }t\geq2.
\end{array}
\right.$\\

Let $\Pi$ be a partition on $V(G\odot H)$ induced by $c$. Now, we will show that the color codes of all vertices are distinct. Let $x$ and $y$ be two vertices of $G\odot H$ such that $c(x)=c(y)$. Since $c((u))\leq l$ and $c((v,b))\geq l+1$, there is no possibility of $x=(u)$ and $y=(v,b)$. So, we only have two possibilities of $x$ and $y$.

\begin{enumerate}
\item $x=(u)$ and $y=(v)$

Since $c(x)=f(u)$ and $c(y)=f(v)$, and $f$ is a locating coloring of $G$, then $c_\Pi(x)\neq c_\Pi(y)$.

\item If $x=(u,a)$ and $y=(v,b)$

The only possibility is $u=v$ and $a,b\in V(H_t)$ for $t\in\{1,2,\ldots,k\}$. Since $c_t$ is a locating coloring of $H_t$, then $a$ and $b$ are differed by $c_t$ which implies $c_{\Pi_t(u)}(x)\neq c_{\Pi_t(u)}(y)$. Therefore, $c_\Pi(x)\neq c_\Pi(y)$.
\end{enumerate}

By all cases above, we obtain that $c$ is a locating coloring of $G\odot H$. Therefore, $\chi_L(G\odot H)\leq \chi_L(G)+\sum_{t=1}^k(\chi_L(H_t+K_1)-1).$ 

Applying Lemmas \ref{teo1bawah} and \ref{teo1atas} above, we obtain the general boundaries of $\chi_L(G\odot H)$ as stated below.

\begin{theorem}\label{teo1}
Let $G$ be a connected graph of order $n\geq2$ and $H$ be an arbitrary graph containing $k$ components $H_1,H_2,\ldots,H_k$. Then $\max\{\chi_L(H_t+K_1)\mid 1\leq t\leq k\} \leq \chi_L(G\odot H)\leq \chi_L(G)+\sum_{t=1}^k(\chi_L(H_t+K_1)-1).$
\end{theorem}

The following two theorems show the existences of graphs $G$ and $H$ which satisfy lower bound and upper bound of Theorem \ref{teo1}, respectively.

\begin{theorem}
There exists a connected graph $G$ of order $n\geq2$ and a graph $H$ containing $k$ components $H_1,H_2,\ldots,H_k$ such that $\chi_L(G\odot H)=\max\{\chi_L(H_t+K_1)\mid 1\leq t\leq k\}$.
\end{theorem}

\noindent \emph{Proof} :

Let $G$ be a path graph with 3 vertices $P_3$ and $H$ be a union of a path graph with 2 vertices $P_2$ and a cycle with 4 vertices $C_4$. Note that, $\chi_L(P_2+K_1)=3$ and $\chi_L(C_4+K_1)=5$. We will show that $\chi_L(G\odot H)=\max\{\chi_L(P_2+K_1),\chi_L(C_4+K_1)\}=5$. By Theorem \ref{teo1}, we only need to show that $\chi_L(G\odot H)\leq\max\{\chi_L(P_2+K_1),\chi_L(C_4+K_1)\}$. Now, we will construct a 5-coloring $c$ in $G\odot H$ such that $c$ is a locating coloring of $G\odot H$.

Let $V(G)=\{u,v,w\}$ with $uv,vw\in E(G)$ and $V(H)=\{a,b,p,q,r,s\}$ with $ab,pq,ps,$ $qr,rs\in E(H)$. We define a 5-coloring $c$ of $G\odot H$ as follows.

$c(y)=\left\{
\begin{array}{ll}
1, & \text{for }y\in\{(v),(u,p),(w,p)\},\\
2, & \text{for }y\in\{(u,q),(v,q),(w,r)\}\cup\{(z,a)|z\in V(G)\},\\
3, & \text{for }y\in\{(w),(u,r),(v,p)\},\\
4, & \text{for }y\in\{(u,s),(v,r),(w,q)\}\cup\{(z,b)|z\in V(G)\},\\
5, & \text{for }y\in\{(u),(v,s),(w,s)\}.
\end{array}
\right.$

\begin{figure}[hbtp]
  \centering
  \includegraphics[width=12cm]{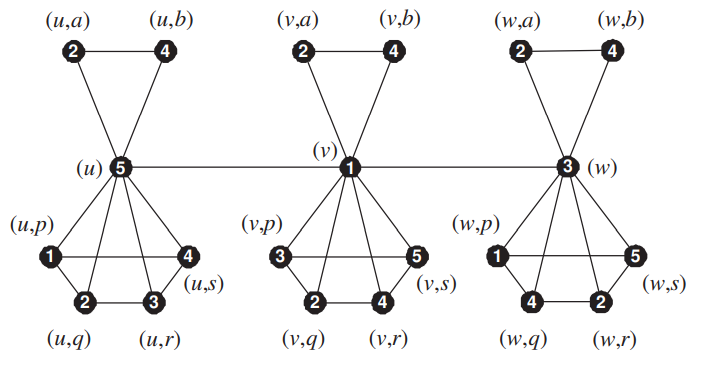}\\
  \caption{A 5-coloring $c$ in $G\odot H$.}
  \label{gbr01}
\end{figure}

Now, we will show that $c$ is a locating coloring of $G\odot H$. Let $\Pi$ be a partition of $V(G\odot H)$ induced by $c$. The color codes of all vertices of $G\odot H$ are as follows.

$\begin{array}{lllll}
c_\Pi((u))=(1,1,1,1,0) & & c_\Pi((v))=(0,1,1,1,1) & & c_\Pi((w))=(1,1,0,1,1)\\
c_\Pi((u,a))=(2,0,2,1,1)& &c_\Pi((v,a))=(1,0,2,1,2)& &c_\Pi((w,a))=(2,0,1,1,2)\\
c_\Pi((u,b))=(2,1,2,0,1)& &c_\Pi((v,b))=(1,1,2,0,2)& &c_\Pi((w,b))=(2,1,1,0,2)\\
c_\Pi((u,p))=(0,1,2,1,1)& &c_\Pi((v,p))=(1,1,0,2,1)& &c_\Pi((w,p))=(0,2,1,1,1)\\
c_\Pi((u,q))=(1,0,1,2,1)& &c_\Pi((v,q))=(1,0,1,1,2)& &c_\Pi((w,q))=(1,1,1,0,2)\\
c_\Pi((u,r))=(2,1,0,1,1)& &c_\Pi((v,r))=(1,1,2,0,1)& &c_\Pi((w,r))=(2,0,1,1,1)\\
c_\Pi((u,s))=(1,2,1,0,1)& &c_\Pi((v,s))=(1,2,1,1,0)& &c_\Pi((w,s))=(1,1,1,2,0)\\
\end{array}
$

Since every two distinct vertices has different color code, $c$ is a locating coloring of $G\odot H$.

\begin{theorem}
There exists a connected graph $G$ of order $n\geq2$ and a graph $H$ containing $k$ components $H_1,H_2,\ldots,H_k$ such that $\chi_L(G\odot H)=\chi_L(G)+\sum_{t=1}^k(\chi_L(H_t+K_1)-1)$.
\end{theorem}

\noindent \emph{Proof} :

Let both $G$ and $H$ be a path graph with two vertices $P_2$. Note that, $\chi_L(P_2)=2$ and $\chi_L(P_2+K_1)=3$. We will show that $\chi_L(G\odot H)=\chi_L(G)+\chi_L(H+K_1)-1=4$. By Theorem \ref{teo1}, we only need to show that $\chi_L(G\odot H)\geq4$.

Suppose that $\chi_L(G\odot H)\leq\chi_L(G)+\chi_L(H+K_1)-2=3$. Since $G\odot H$ contains $C_3$, $G\odot H$ needs at least 3 colors where three vertices in $C_3$ have distinct colors. So, it follows that $\chi_L(G\odot H)=3$. Since there are two disjoint $C_3$ in $G\odot H$, two vertices with color $i\in\{1,2,3\}$ have the same color code, a contradiction. Therefore, $\chi_L(G\odot H)\geq 4$. 

In theorem below, we also give an existences of graphs $G$ and $H$ which do not satisfy both lower bound and upper bound of Theorem \ref{teo1}.

\begin{theorem}\label{teo2}
There exists a connected graph $G$ of order $n\geq2$ and a graph $H$ containing $k$ components $H_1,H_2,\ldots,H_k$ such that $\chi_L(G\odot H)=l$ where $\max\{\chi_L(H_t+K_1)\mid 1\leq t\leq k\} \leq l\leq \chi_L(G)+\sum_{t=1}^k(\chi_L(H_t+K_1)-1).$.
\end{theorem}

\noindent \emph{Proof} :

Let $H$ be a complement of complete graph with $k\geq 2$ vertices $\overline{K_k}$ and $G$ be a connected graph with $n\leq k+1$ vertices. In this case, $H$ contains $k$ components of one vertex $P_1$ and $\chi_L(P_1+K_1)=2$. Now, we will show that $\chi_L(G\odot H)=k+1$. Note that, $\max\{\chi_L(H_t+K_1)\mid 1\leq t\leq k\}=2<k+1<\chi_L(G)+k=\chi_L(G)+\sum_{t=1}^k(\chi_L(H_t+K_1)-1)$.

Since $G\odot H$ contains a vertex which is adjacent to $k$ endpoints, by Corollary \ref{ColA}, $\chi_L(G\odot H)\geq k+1$.

Let $V(G)=\{u_1,u_2,\ldots,u_n\}$ and $V(H)=\{v_1,v_2,\ldots,v_k\}$. Now, we define $c:V(G\odot H)\rightarrow \{1,2,\ldots,k+1\}$ as follows.

$c((u_i))=i$ for $1\leq i\leq n$, and

$c((u_i,v_j))=\left\{
\begin{array}{ll}
j, & \text{for }1\leq i\leq n\text{ and }1\leq j\leq k\text{ and }i\neq j,\\
k+1, & \text{for }1\leq i\leq n\text{ and }1\leq j\leq k\text{ and }i=j.
\end{array}
\right.
$

We will show that $c$ is a locating coloring of $G\odot H$. Let $\Pi=\{C_1,C_2,\ldots,C_{k+1}\}$ be a partition of $V(G\odot H)$ induced by $c$. Let $x$ and $y$ be two vertices of $G\odot H$ having same color. We have two cases of $x$ and $y$.

\begin{enumerate}
\item $x=(u_i)$ and $y=(u_p,v_j)$ with $p\neq i$

Since $x$ is adjacent to vertex with color $k+1$ and $c((u_i))\neq k+1$ for $1\leq i\leq n$, we obtain that $d(x,C_{k+1})\neq d(y,C_{k+1})$. It follows that $c_\Pi(x)\neq c_\Pi(y)$.

\item $x=(u_i,v_j)$ and $y=(u_p,v_q)$ with $p\neq i$ and $j,q\in \{1,2,\ldots,k\}$

Since $d(x,C_i)=1\neq2=d(y,C_i)$, we obtain that $c_\Pi(x)\neq c_\Pi(y)$.
\end{enumerate}

Therefore, $c$ is a locating coloring of $G\odot H$.

\section{The  Locating-Chromatic Number of $T_n\odot \overline{K}_m$}

In this section, we consider the locating-chromatic number of $T_n\odot \overline{K_m}$ where $T_n$ is a tree graph with $n$ vertices and $\overline{K_m}$ is complement of complete graph with $m$ vertices. On the other hand, $\overline{K_m}$ is a graph without edges. First, we give the boundaries of locating-chromatic number of $T_n\odot \overline{K}_m$ as a special case of Theorem \ref{teo1} and a direct consequences of Corollary \ref{ColA}.

\begin{theorem}\label{teo3}
For $m\geq 1$ and $n\geq2$, let $T_n$ be a tree with $n$ vertices and $\overline{K_m}$ be a complement of complete graph with $m$ vertices. Then $m+1\leq \chi_L(T_n\odot \overline{K}_m)\leq \chi_L(T_n)+m$.
\end{theorem}

The existence of lower bound for Theorem \ref{teo3} can be seen in the next theorem. The proof of this theorem is similar with the proof of Theorem \ref{teo2}.

\begin{theorem}
For $m,n\geq2$, let $T_n$ be a tree with $n$ vertices and $\overline{K_m}$ be a complement of complete graph with $m$ vertices. If $1\leq n \leq m+1$, then $ \chi_L(T_n\odot \overline{K}_m)=m+1$.
\end{theorem}

\begin{figure}[hbtp]
  \centering
  \includegraphics[width=14cm]{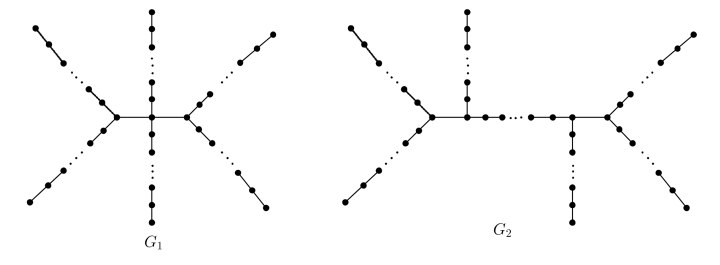}\\
  \caption{Graphs $G_1$ and $G_2$}
  \label{gbr02}
\end{figure}

For an existence of upper bound of Theorem \ref{teo3} and an existence of graph $T_n$ with $\chi_L(T_n\odot \overline{K_m})$ is not equal to both upper and lower bound of Theorem \ref{teo3}, we consider $T_n\odot \overline{K_1}$. Note that, $T_n\odot \overline{K_1}$ is isomorph to $T_n\odot K_1$. Generally, for graph $G$ with $n\geq2$ vertices, $\chi_L(G)\geq2$. Furthermore, $\chi_L(G)=2$ if and only if $|V(G)|=2$. In \cite{Baskoro2013}, Baskoro \emph{et al.} have characterized all trees with locating-chromatic number $3$, which is subgraph of graphs $G_1$ or $G_2$ (Figure \ref{gbr02}).  Motivated by this, since $T_n\odot \overline{K_1}$ is a class of tree, in this paper we give a condition of $T_n$ with $\chi_L(T_n)=3$ such that the locating-chromatic number of $T_n\odot \overline{K}_1$ is equal to 3 or 4.

\begin{figure}[hbtp]
  \centering
  \includegraphics[scale=.6]{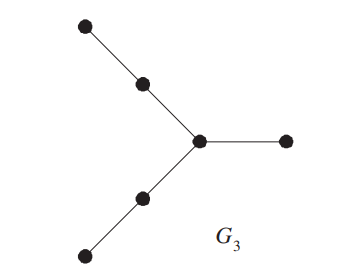}\\
  \caption{Graph $G_3$}
  \label{gbr03}
\end{figure}

\begin{theorem}
For $n\geq 2$, let $T_n$ be a tree with $n$ vertices such that $\chi_L(T_n)=3$. If $G_3$ is a graph as stated in Figure \ref{gbr03}, then
\begin{equation*}
\chi_L(T_n\odot \overline{K_1})=\left\{
\begin{array}{ll}
3, & \text{if } T_n\text{ is a subgraph of }P_6\text{ or }G_3,\\
4, & \text{otherwise}.
\end{array}
\right.
\end{equation*}

For $n\geq 2$, let $T_n$ be a tree with $n$ vertices such that $\chi_L(T_n)=3$. If $T_n$ is a subgraph of $P_6$ or $G_3$ (Figure \ref{gbr03}), then $\chi_L(T_n\odot \overline{K_1})=3.$
\end{theorem}

\noindent \emph{Proof} :

We distinguish two cases.

\noindent \textbf{Case 1.} $T_n$ is a subgraph of $P_6$ or $G_3$.

Then $T_n\odot \overline{K_1}$ is a graph which is a subgraph of $G_1$ or $G_2$. By \cite{Baskoro2013}, $\chi_L(T_n\odot \overline{K_1})=3$.

\noindent \textbf{Case 2.} $T_n$ is not a subgraph of $P_6$ or $G_3$.

Then $T_n\odot \overline{K_1}$ is not a subgraph of $G_1$ or $G_2$. So, $\chi_L(T_n\odot \overline{K_1})\geq4$. Since from Theorem \ref{teo3}, $\chi_L(T_n\odot \overline{K_1})\leq \chi_L(T_n)+1=4$. Therefore, $\chi_L(T_n\odot \overline{K_1})=4$. 

We also give an additional existence of tree graph $T_n$, namely star graph such that $\chi_L(T_n\odot \overline{K_m})$ is not equal to both lower and upper bound in Theorem \ref{teo3} for some $m\geq 1$. Note that, a star graph with $n$ vertices $S_n$ satisfies $\chi_L(S_n)=n$.

\begin{theorem}
For $n\geq4$, let $S_n$ be a star with $n$ vertices. Then $\chi_L(S_n\odot \overline{K_1})=\lceil\sqrt{n}\rceil+1.$
\end{theorem}

\noindent \emph{Proof} :

Let $V(S_n\odot \overline{K_1})=\{x,y\}\cup\{x_i,y_i\mid 1\leq i\leq n-1\}$ and $E(S_n\odot \overline{K_1})=\{xy\}\cup\{xx_i,x_iy_i\mid1\leq i\leq n-1\}$.

Suppose that $\chi_L(S_n\odot K_1)\leq\lceil\sqrt{n}\rceil=k$. Let $c$ be a locating coloring of $S_n\odot\overline{K_1}$ with $k$ colors and $\Pi$ be a partition of $V(S_n\odot\overline{K_1})$ induced by $c$. Without loss of generality, let $c(x)=1$. For vertices of $A=\{z\in V(S_n\odot\overline{K_1})\mid xz\in E(S_n\odot\overline{K_1})\}$, we have $k-1$ possibility of colors. Since $|A|=n$, there exist at least $k$ vertices of $A$ having same color. Without loss of generality, let $A_1=\{x_1,x_2,\ldots,x_k\}$ be the vertices of $A$ having same color. Note that, for two distinct vertices $x_i$ and $x_j$ of $A_1$ with $1\leq i<j\leq k$, $d(x_i,z)=d(x_j,z)$ for $z\in V(S_n\odot\overline{K_1})\setminus\{x_i,x_j,y_i,y_j\}$. Since $d(x_i,y_i)\neq d(x_j,y_i)$, we have $c(y_{i})\neq c(y_j)$. Since there are $k-1$ possible colors for $\{y_1,y_2,\ldots,y_k\}$, there exist two different vertices $y_i$ and $y_j$ with $i,j\in\{1,2,\ldots,k\}$ and $i\neq j$ such that $c(y_i)=c(y_j)$. It follows that $c_\Pi(x_i)=c_\Pi(x_j)$, a contradiction.

Now, we will show that $\chi_L(S_n\odot \overline{K_1})\leq\lceil\sqrt{n}\rceil+1=l$. For $1\leq t\leq \lceil\frac{n}{l-1}\rceil$, $(t-1)(l-1)< i\leq t(l-1)$, and $1\leq j\leq l-1$, we define an $l$-coloring $c:V(S_n\odot \overline{K}_1)\rightarrow \{1,2,\ldots,l\}$ as follows.

$c(x)=1$, $c(y)=l$, $c(x_i)=t+1$, and

$c(y_{(t-1)(l-1)+j})=\left \{\begin{array}{ll}
          l-j+1, & \textrm{if }l-j>t, \\
          l-j, & \textrm{otherwise.}
        \end{array} \right.$

We will show that $c$ is a locating coloring of $S_n\odot \overline{K_1}$.  Let $\Pi$ be a partition of $V(S_n\odot\overline{K_1})$ induced by $c$. Let $u$ and $v$ be two distinct vertices of $S_n\odot \overline{K_1}$ such that $c(u)=c(y)$. We distinguish six cases.

\noindent \textbf{Case 1}. $u=x$ and $v=y_i$ for $i\in\{1,2,\ldots,n-1\}$

Since $v$ is only adjacent to $x_i$ and $u$ is adjacent to vertices with colors $2,3,\ldots,l$, we obtain that $c_\Pi(u)\neq c_\Pi(v)$.

\noindent \textbf{Case 2}. $u=y$ and $v=y_i$ for $i\in\{1,2,\ldots,n-1\}$

Since $u$ is only adjacent to $x$ where $c(x)=1$ and $v$ is adjacent to $x_i$ where $c(x_i)\in\{2,3,\ldots,l-1\}$, we obtain that $c_\Pi(u)\neq c_\Pi(v)$.

\noindent \textbf{Case 3}. $u=y$ and $v=x_i$ for $i\in\{1,2,\ldots,n-1\}$

Since $u$ is only adjacent to $x$ and $v$ is adjacent to $x$ and $y_i$, we obtain that $c_\Pi(u)\neq c_\Pi(v)$.

\noindent \textbf{Case 4}. $u=x_i$ and $v=x_j$ for $i,j\in\{1,2,\ldots,n-1\}$ and $i\neq j$

The vertex $u$ is adjacent to $x$ and $y_i$, and the vertex $v$ is adjacent to $x$ and $y_j$. Since $c(y_i)\neq c(y_j)$, we obtain that $c_\Pi(u)\neq c_\Pi(v)$.

\noindent \textbf{Case 5}. $u=y_i$ and $v=y_j$ for $i,j\in\{1,2,\ldots,n-1\}$ and $i\neq j$

Since $u$ and $v$ are only adjacent to $x_i$ and $x_j$, respectively, and $c(x_i)\neq c(x_j)$, we obtain that we obtain that $c_\Pi(u)\neq c_\Pi(v)$.

\noindent \textbf{Case 6}. $u=x_i$ and $v=y_j$ for $i,j\in\{1,2,\ldots,n-1\}$

Note that, $i\neq j$. So, $v$ is only adjacent to $x_j$ which is not $u$. By the definition of $c$ above, $c(x_j)\neq c(v)$. Since $u$ is adjacent to $y_i$ which is not $v$, and $x$, we obtain that $c_\Pi(u)\neq c_\Pi(v)$.

By cases above, we obtain that $c$ is a locating coloring of $S_n\odot\overline{K_1}$. 

%

\end{document}